\documentclass{article}

\usepackage{amsmath, amsthm, amsfonts}

\theoremstyle{definition}

\theoremstyle{theorem}

\theoremstyle{definition}

\theoremstyle{remark}

\usepackage{graphicx}
\usepackage{subfigure}          

\usepackage{placeins}

\usepackage{tabularx}

\usepackage{enumitem}

\usepackage{lmodern,textcomp}

\usepackage{algorithm}
\usepackage{algpseudocode}

\usepackage{enumitem}
\newlist{Aufz}{enumerate}{10}
\setlist[Aufz]{label*=\arabic*.}

\usepackage{xspace} 
\newcommand\largeparbreak{\par\bigskip}

\newcommand{\R}{\mathbb{R}\xspace}
\newcommand{\N}{\mathbb{N}\xspace}

\newcommand{\cO}{\mathcal{O}\xspace}

\newcommand{\bA}{\textbf{A}\xspace}

\newcommand{\bB}{\textbf{B}\xspace}

\newcommand{\tbY}{\tilde{\textbf{Y}}\xspace}

\newcommand{\bU}{\textbf{U}\xspace}
\newcommand{\bV}{\textbf{V}\xspace}
\newcommand{\bZ}{\textbf{Z}\xspace}
\newcommand{\bS}{\textbf{S}\xspace}

\newcommand{\bD}{\textbf{D}\xspace}

\newcommand{\bx}{\textbf{x}\xspace}

\newcommand{\by}{\textbf{y}\xspace}

\newcommand{\hbY}{\hat{\textbf{Y}}\xspace}

\newcommand{\tbV}{\tilde{\textbf{V}}\xspace}

\newcommand{\bM}{\textbf{M}\xspace}

\newcommand{\bC}{\textbf{C}\xspace}
\newcommand{\bX}{\textbf{X}\xspace}
\newcommand{\bW}{\textbf{W}\xspace}
\newcommand{\bY}{\textbf{Y}\xspace}

\newcommand{\bI}{\textbf{I}\xspace}

\newcommand{\bO}{\textbf{0}\xspace}

\title{A time-optimal algorithm for solving (block-)tridiagonal linear systems of dimension $N$ on a distributed computer of $N$ nodes}

\author{Martin Neuenhofen}

\date{\today}

\begin{document}

\maketitle

\begin{abstract}
	We are concerned with the fastest possible direct numerical solution algorithm for a thin-banded or tridiagonal linear system of dimension $N$ on a distributed computing network of $N$ nodes that is connected in a binary communication tree. Our research is driven by the need for faster ways of numerically solving discretized systems of coupled one-dimensional black-box boundary-value problems.
	
	Our paper presents two major results: First, we provide an algorithm that achieves the optimal parallel time complexity for solving a tridiagonal linear system and thin-banded linear systems. Second, we prove that it is impossible to improve the time complexity of this method by any polynomial degree.
	
	To solve a system of dimension $m\cdot N$ and bandwidth $m \in \Omega(N^{1/6})$ on $2 \cdot N-1$ computing nodes, our method needs time complexity $\cO(\log(N)^2 \cdot m^3)$.
\end{abstract}


\section{Introduction}

\paragraph{Motivation of the problem}
Many computational engineering tasks deal with the solution of (systems) of one-dimensional differential-algebraic boundary-value problems \cite{Russell1972}. Examples are numerical simulations of the following physical phenomena: the deformation of a clamped beam, the dynamic pressure in a gas-pipe, the trajectory of a missile, and constrained optimal control problems.

Using a numerical discretization method, a large-dimensional equation systems results that is typically solved via Newton's method. The arising linear systems are thin-banded and of very large dimension.

Often, the Newton system results from a minimization principle, either of an objective or by a natural model that aims for minimization of potential energy. In these cases the arising linear systems are not only thin-banded and of very large dimension, but they are also symmetric positive definite, which is clearly desirable for reasons of numerical stability.

Since thin-banded, the system matrix can be interpreted as block-tridiagonal, where the block-size is identical to the band-width. Thus, for ease of presentation in the following we present an algorithm for (block-) tridiagonal systems where the block-size is $m\ll N$.

\paragraph{Motivation of the problem statement and computing model}
In consequence of the emergence of massively parallel computing systems, nowadays the problem in numerical computing is typically not to solve a given mathematical problem, but rather to solve it on a given computing system while exploiting its resources in an optimal way. Especially for parallel computing systems it is difficult to spread the computational task in a way that enables the full use of the computing system's capacity.

As for our case of solving thin-banded linear systems, it 
is well-known that for a bandwidth bounded by $m \in \cO(1)$ the optimal time complexity for solving a banded linear system on a serial computer is
\begin{align*}
	\cO(N)\,,
\end{align*}
as can be achieved by use of Gaussian elimination. This is equivalent in meaning to that solving the problem with a serial machine is literally as expensive as the following two: \textit{reading} the problem, or \textit{writing} the problem's solution into the memory.

So at first glance it seems like nothing can be done to improve: The time to communicate the problem to a solver would already predominate the time that the solver actually needed. So where is the point in trying to make the solver faster?

The point or answer is that problems do not need to be communicated. One could assemble in parallel the rows of a linear equation system in a distributed way on the memory of independent computing nodes that are connected via a network. Using the algorithm that we propose, all the computing nodes can solve the one big linear system but each of them does only \textit{read a tiny part} of the problem and only \textit{writes a tiny part} of the solution vector.

\paragraph{Literature review}
We consider the problem of solving a thin-banded linear system as a generalization of solving tridiagonal linear systems in parallel, which is why our literature review refers to parallel tridiagonal solvers. Such solvers can be applied by interpreting the thin-banded linear system as block-tridiagonal system of dense blocks.

There are several popular methods for the solution of tridiagonal linear systems. These have in common that they use a concept described as \textit{parallel factorizations} \cite{Amodio92parallelfactorizations,MATTOR19951769}: The matrix is multiplied from the left with a block-diagonal matrix that decouples a portion of the unknowns through a n interface system. The reduced system is solved. The reduced solution is distributed to all processors so that they can compute in parallel the formerly removed unknowns.

According to \cite{austin-2004-linesolver}, the first parallel tridiagonal solver is called \textit{cyclic reduction} and was presented 1965 in \cite{Hockney:1965:FDS:321250.321259}. The \textit{recursive doubling algorithm} was introduced 1973 in \cite{Stone:1973:EPA:321738.321741}. In both of these algorithms each processor holds one row of the system. Cyclic reduction works by successively expressing the odd variables of the solution vector in terms of the even. This can be done until finally there is a system of one variable that is solved. An implementation is provided in \cite{DBLP:journals/siamsc/BrownFJ00}.

\textit{Wang's method} \cite{Wang:1981:PMT:355945.355947}, introduced 1981, is a parallel algorithm where the order of the number of processors is smaller than the order of the number of rows of the system to be solved. This algorithm has been proven to be numerically stable \cite{Yalamov:1999:SPA:300088.300141}. The idea of this method is the assembly of an interface problem whose solution can be used directly to solve the remaining variables in a backward-substitution step.

In 1991 Bondeli introduced a divide and conquer algorithm for tridiagonal system \cite{Bondeli:1991:PDC:1746085.1746141}. The idea of this algorithm is to solve a block-diagonal approximation of the system in parallel, where each processor holds a diagonal-block. It results a reduced interface system that is solved by cyclic reduction, cf. \cite{austin-2004-linesolver}.

\paragraph{Organization of the paper}
In the remainder of this section we describe the parallel computing system that we use for our algorithm. We then describe the mathematical problem that the algorithm solves. It is important that the problem is given to the algorithm in a special way. In particular, the problem data must be provided in distributed memory before the algorithm is called. This is important because it would take to much time to move the data from a central storage to the distributed memory.

In Section~2 we present the algorithm. We start with an implementation and afterwards show how this algorithm can be derived from familiar matrix algorithms. At the end of the section we analyse the time complexity of the algorithm and remark on optimality of this complexity result.

Section 3 we give lower bounds on the time complexity for solving tridiagonal linear systems on a parallel distributed memory machine. We show that our algorithm is able to yield optimal time complexity.

Eventually we draw conclusions in Section~4.

\paragraph{Computing system}

We need to describe the computing system before we describe the problem statement because otherwise we cannot describe where we presume the problem data to be placed. Above we described why this is crucial: We have to make sure that the problem is provided in the right way because moving the problem data around would cost too much time.

\begin{figure}
\centering
\includegraphics[width=1\linewidth]{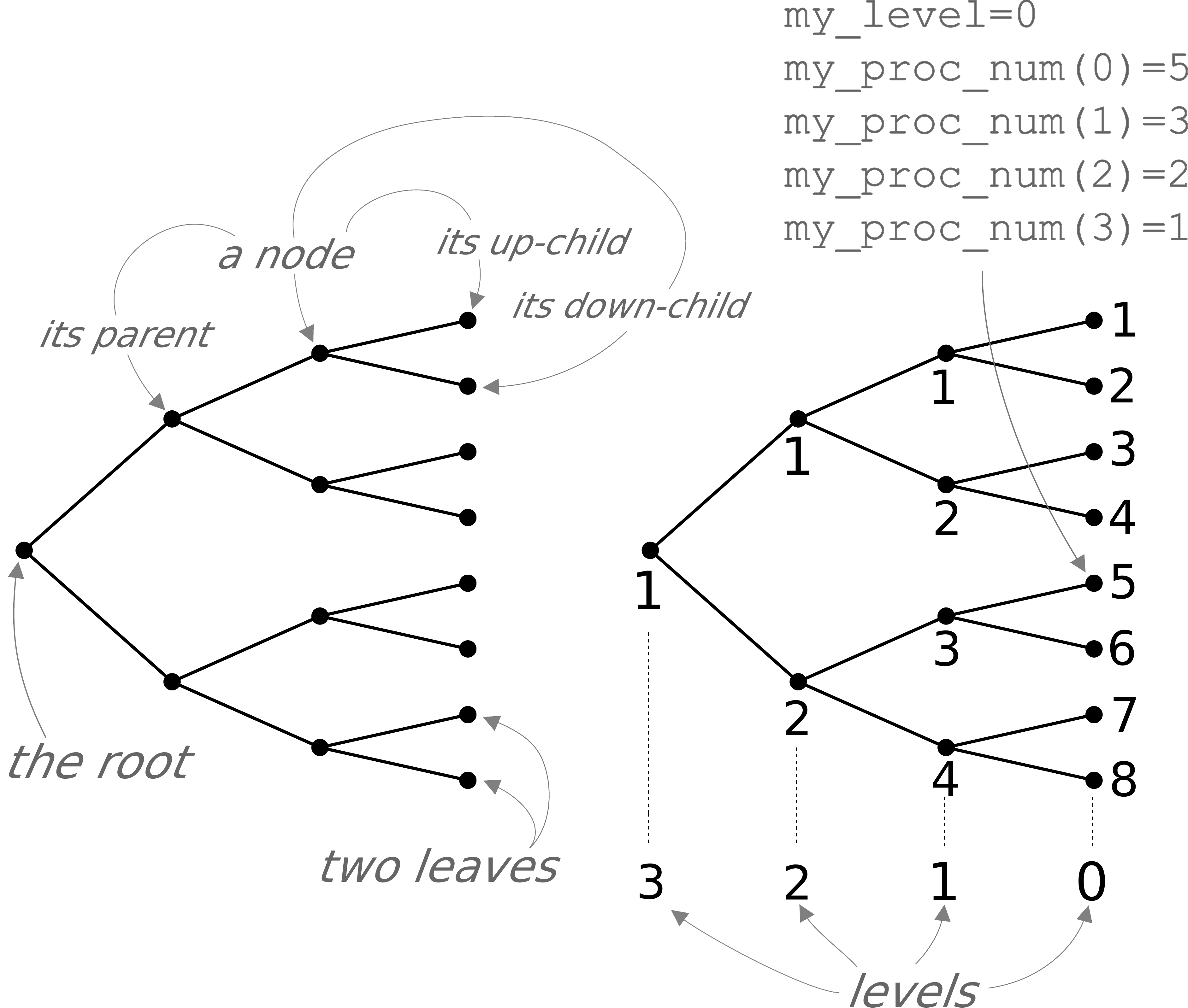}
\caption{Structure of the computing system: $2 \cdot N-1$ nodes are connected in a two-tree communication network. Each node can send messages to its parent and children. The nodes are classified in levels. Further, each node per level is given a processor number. Each node knows its level and the processor number of itself and its parents.}
\label{fig:ComputingSystem}
\end{figure}

For a system of dimension $N \cdot m$ with a bandwidth $m$ we consider a computing system of $2 \cdot N-1$ computing nodes that have each their own memory and that each run the solution algorithm in parallel. The nodes are connected via a network of cables. Nodes can send data to other nodes and nodes can wait to receive data from others. Figure~\ref{fig:ComputingSystem} illustrates the computing system. The black circles symbolize the computing nodes and the black lines illustrate their connections via a cable network.

For the network we require a special structure: For our algorithm we need a two-tree network. In a two-tree, also called dual tree, each node is connected to three other nodes: a \textit{parent}, an \textit{up-child} and a \textit{down-child}, cf. in the figure. Exceptions are: There exists one node called \textit{root}. This node does not have a parent. Further, there exist $N$ nodes that are called \textit{leaves}. A leaf does only have a parent, but it does neither have an \textit{up-child} nor a \textit{down-child}.

The right part of the figure assigns the nodes with numbers. Each node has a \textit{processor number} and a \textit{level}. The levels are defined recursively: Each leaf has a level of zero, and the parent of each node has a level that is by one larger than the level of the node itself. The processor number is a value that is given by counting from 1 from the uppermost node to the lowermost node of each level.

Each node holds identifying variables like a passport: The variable \texttt{my\_level} gives the level of this node. The variable \texttt{my\_proc\_num} is a list. The following values are well-defined
\begin{align*}
	\texttt{my\_proc\_num}(\ell) 	\quad \text{ for }\ell \in \lbrace \texttt{my\_level},...d\rbrace\,,
\end{align*}
where $d$ is the level of the root. The value $\texttt{my\_proc\_num}(\ell)$ gives the processor number of the node on level $\ell$ through which a signal from the root would have to travel in order to reach this node. Figure~\ref{fig:ComputingSystem} gives an example for this: To reach the node on level $0$ with processor number $5$, a signal from the root would have to traverse through node $2$ of level $2$ and node $3$ of level $1$.

We presume the nodes as identical serial computing units. As is common, we presume that basic operations plus, minus, times, divide, and copy of scalar values require each a fixed amount of time on a respective node.

For communications over the network, we use the following message passing interface of six commands:
\begin{itemize}
	\item \texttt{send\_to\_up\_child}($\bM$). If this node has an up-child then it sends a matrix by copy $\bM$ to up-child and waits until up-child received it.
	\item \texttt{receive\_from\_up\_child}($\bM$). If this node has a down-child then it waits until it receives by copy a matrix from up-child. This node stores the matrix in its variable $\bM$ and the copy for the transmission is destroyed.
	\item \texttt{send\_to\_down\_child}($\bM$); analogous to above, but the data is sent to the down-child of this node.
	\item \texttt{receive\_from\_down\_child}($\bM$); analogous to above, but the data is received from the down-child of this node.
	\item \texttt{send\_to\_parent}($\bM$); analogous to above, but the data is sent to the parent of this node.
	\item \texttt{receive\_from\_parent}($\bM$); analogous to above, but the data is received from the parent of this node.
\end{itemize}
For each communication we assume a time complexity of the number of elements of $\bM$ plus a constant amount of time $c^\text{Lat}_N$ that is due to latency. The latency accounts for the phenomenon that information travels through the cable at speed of light, so it takes a while until the beginning of a message has moved through the cable.

\paragraph{Problem statement}
We consider the numerical solution of a banded linear system
\begin{align}
\underline{\bA} \cdot \underline{\bX} = \underline{\bY}\label{eqn:LinearSystem}
\end{align}
where $\underline{\bA} \in \R^{(N\cdot m) \times (N \cdot m)}$ has bandwidth $b\leq m$, and $\underline{\bY} \in \R^{(N \times m) \times k}$ is a dense matrix of $k$ right-hand sides. The task is to find numerical values for $\underline{\bX} \in \R^{(N \times m) \times k}$. We assume $N \in 2^\N$.
\largeparbreak

As formerly discussed, the time for writing the solution into memory in a sequential way would already exceed the time that is actually needed to solve the system. This is why in the following we describe very precisely in which form the data $\underline{\bA}$, $\underline{\bY}$ must be provided to our computing system.

\begin{figure}
\centering
\includegraphics[width=1\linewidth]{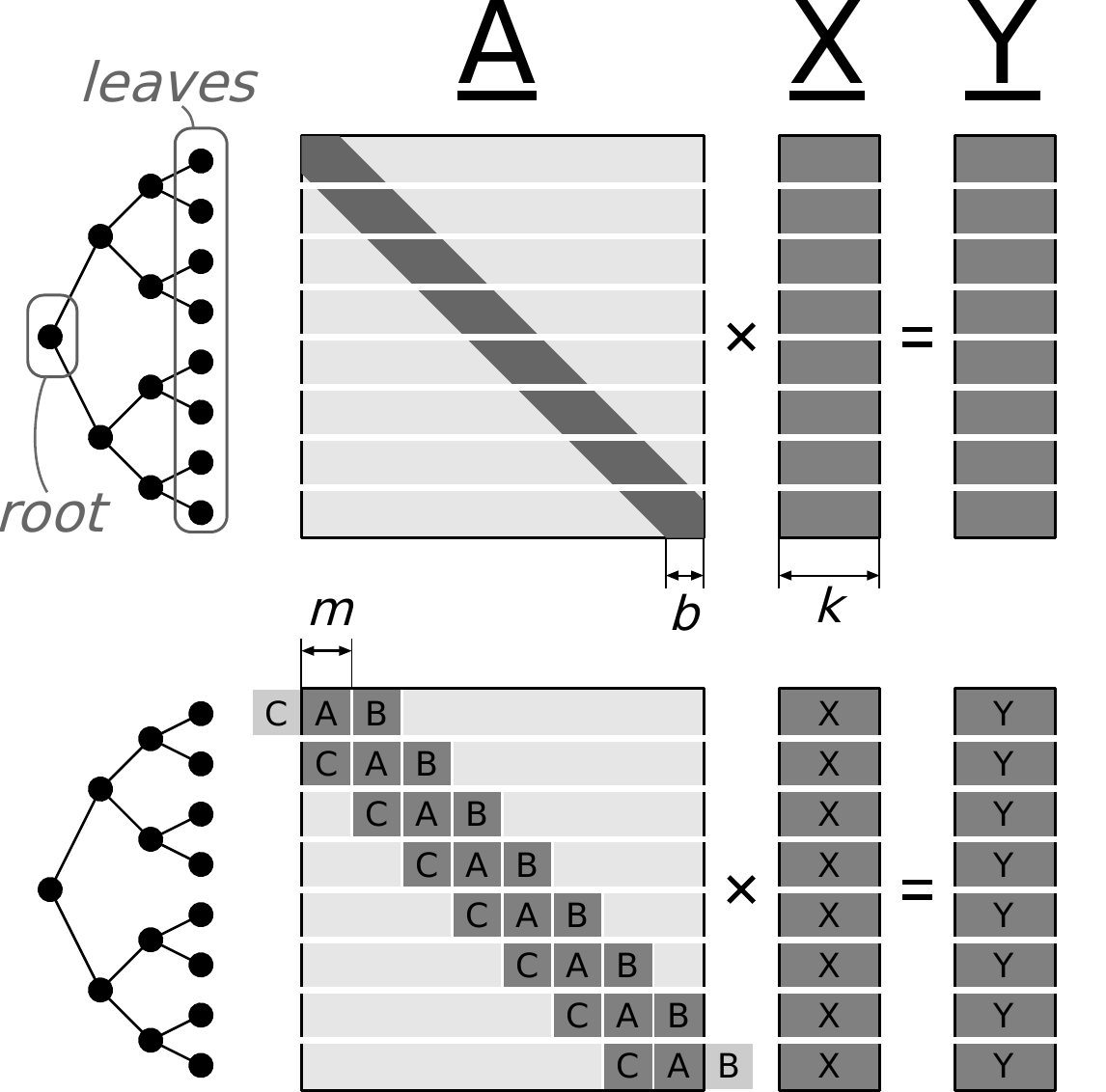}
\caption{Distributed storage of the linear system. At the top: The linear system of bandwidth $b \leq m$ is chunked into row-blocks of size $m$. At the bottom: Each row-block is stored in one leaf using five matrices.}
\label{fig:Block_Matrix_Storage}
\end{figure}

The system matrix, the right-hand sides, and the solution vectors are stored in a separated way in the leaves of our two-tree. Figure~\ref{fig:Block_Matrix_Storage} illustrates the situation for $N=8$. Each leaf holds five matrices in its private storage: $\bA,\bB,\bC \in \R^{m \times m}$ and $\bX,\bY \in \R^{m \times k}$. Matrices of two distinct leaves can have totally different values.
Comparing the upper and lower part of the Figure~\ref{fig:Block_Matrix_Storage}, we find that the original matrices $\underline{\bA}$, $\underline{\bX}$, $\underline{\bY}$ can be composed of the matrices $\bA,\bB,\bC,\bX,\bY$ of all leaves. There are two matrices that fall out of the pattern: $\bC$ in the uppermost leaf and $\bB$ in the lowermost leaf. We require that these matrices are zero-matrices.

\section{The algorithm}

This section is organized as follows. We first present the algorithm as a code that could be directly used for implementation in a programming language such as MPI with Cpp or Fortran. Then we sketch our derivation of the algorithm, that arised from applying the devide-and-conquer paradigm on the SPIKE algorithm due to Sameh \cite{Polizzi2007}. We explain how our algorithm operates and give an example to illustrate the algorithmic steps. Finally, we analyse the parallel time complexity.

\paragraph{The algorithm}
The following algorithm is launched on all nodes of the computing system at the same time with their respective local data.
\begin{algorithmic}[1]
		\Procedure{ParallelSolver}{$\bA,\bB,\bC,\bY,N,d$}
		\State \textit{// matrices: }$\bV^{\lbrace j\rbrace}\in \R^{m \times m}$ for $j=1,...,d$; $\bV_\text{up},\bV_\text{down} \in \R^{m \times m}$
		\State \textit{// matrices (cont. 1): }$\bZ^{\lbrace j\rbrace}_\text{V,up},\bZ_\text{V}^{\lbrace j \rbrace},\bZ^{\lbrace j \rbrace}_\text{V,down} \in \R^{m \times m}$ for $j=1,...,d$
		\State \textit{// matrices (cont. 2): } $\bZ_\text{X,up},\bZ_\text{X},\bZ_\text{X,down} \in \R^{m \times k}$
		\If{$\texttt{my\_level}==0$}
			\State \textit{// - - - write wings}
			\State $j_B := \texttt{my\_proc\_num(my\_level)}$\ ; \quad $k:=N/2$
			\For{$j=d\ :\ -1\ :\ 1$}
				\If{$j_B\geq k$}
					\State $j_B := j_B - k$
				\EndIf
				\If{$j_B==0$}
					\State $\bV^{\lbrace j\rbrace} := \bV^{\lbrace j\rbrace} + \bB$\ ; \quad \textbf{break for-loop}
				\EndIf
				\State $k:=k/2$
			\EndFor
			\State $j_C := \texttt{my\_proc\_num(my\_level)}-1$\ ; \quad $k:=N/2$
			\For{$j=d\ :\ -1\ :\ 1$}
				\If{$j_C\geq k$}
				\State $j_C := j_C - k$
				\EndIf
				\If{$j_C==0$}
				\State $\bV^{\lbrace j\rbrace} := \bV^{\lbrace j\rbrace} + \bC$\ ; \quad \textbf{break for-loop}
				\EndIf
				\State $k:=k/2$
			\EndFor
			\State \textit{// - - - block-diagonal inversion}
			\State $[\bV^{\lbrace 1\rbrace},...,\bV^{\lbrace d\rbrace},\bX] := \bA \backslash [\bV^{\lbrace 1\rbrace},...,\bV^{\lbrace d\rbrace},\bY]$
		\EndIf
		\For{$\ell=1\ :\ 1\ :\ d$}
			\If{$\texttt{my\_level}==0$}
				\State \texttt{send\_to\_parent}($\bV^{\lbrace\ell\rbrace},...,\bV^{\lbrace d \rbrace}\,,\,\bX$)
			\Else
				\State \texttt{receive\_from\_up\_child}($\bV^{\lbrace\ell\rbrace}_\text{up},...,\bV^{\lbrace d\rbrace}_\text{up}\,,\,\bX_\text{up}$)
				\State \texttt{receive\_from\_down\_child}($\bV^{\lbrace\ell\rbrace}_\text{down},...,\bV^{\lbrace d\rbrace}_\text{down}\,,\,\bX_\text{down}$)
				\If{$\texttt{my\_level}<\ell$}
					\If{$\texttt{my\_proc\_num}(\ell-1)$ is odd}
						\State \texttt{send\_to\_parent}($\bV^{\lbrace\ell\rbrace}_\text{down},...,\bV^{\lbrace d\rbrace}_\text{down}\,,\,\bX_\text{down}$)
					\Else
						\State \texttt{send\_to\_parent}($\bV^{\lbrace\ell\rbrace}_\text{up},...,\bV^{\lbrace d\rbrace}_\text{up}\,,\,\bX_\text{up}$)
					\EndIf
				\EndIf
			\EndIf
			\State \textit{// above: nodes of level $\ell$ receive $\bV^{\lbrace \ell,...,d \rbrace }_\text{up},\bV^{\lbrace \ell,...,d \rbrace }_\text{down},\bX_\text{up},\bX_\text{down}$}
			\If{$\texttt{my\_level}==\ell$}
				\State $\bS := \begin{bmatrix}
					\bI_{m \times m} & \bV_\text{up}^{\lbrace \ell \rbrace} \\
					\bV_\text{down}^{\lbrace \ell \rbrace} & \bI_{m \times m}
				\end{bmatrix}$
				\State \mbox{$\begin{bmatrix}
				\bZ_\text{V,up}^{\lbrace \ell+1 \rbrace},...,\bZ_\text{V,up}^{\lbrace d \rbrace} & \bZ_\text{X,up}\\
				\bZ_\text{V,down}^{\lbrace \ell+1 \rbrace},...,\bZ_\text{V,down}^{\lbrace d \rbrace} & \bZ_\text{X,down}
				\end{bmatrix} := \bS \, \backslash\, \begin{bmatrix}
				\bV_\text{up}^{\lbrace \ell+1 \rbrace},...,\bV_\text{up}^{\lbrace d \rbrace} & \bX_\text{up}\\
				\bV_\text{down}^{\lbrace \ell+1 \rbrace},...,\bV_\text{down}^{\lbrace d \rbrace} & \bX_\text{down}
				\end{bmatrix}$}
				\State \texttt{send\_to\_up\_child}($\bZ_\text{V,down}^{\lbrace \ell+1 \rbrace},...,\bZ_\text{V,down}^{\lbrace d \rbrace}\,,\,\bZ_\text{X,down}$)
				\State \texttt{send\_to\_down\_child}($\bZ_\text{V,up}^{\lbrace \ell+1 \rbrace},...,\bZ_\text{V,up}^{\lbrace d \rbrace}\,,\,\bZ_\text{X,up}$)
			\EndIf
			\If{$\texttt{my\_level}<\ell$}
				\State \texttt{receive\_from\_parent}($\bZ^{\lbrace \ell+1 \rbrace}_\text{V},...,\bZ^{\lbrace d \rbrace}_\text{V}\,,\,\bZ_\text{X}$)
				\If{$\texttt{my\_level}>0$}
					\State \texttt{send\_to\_up\_child}($\bZ^{\lbrace \ell+1 \rbrace}_\text{V},...,\bZ^{\lbrace d \rbrace}_\text{V}\,,\,\bZ_\text{X}$)
					\State \texttt{send\_to\_down\_child}($\bZ^{\lbrace \ell+1 \rbrace}_\text{V},...,\bZ^{\lbrace d \rbrace}_\text{V}\,,\,\bZ_\text{X}$)
				\EndIf
				\If{$\texttt{my\_level}==0$}
					\State \mbox{$[\, \bV^{\lbrace \ell+1 \rbrace},...,\bV^{\lbrace d \rbrace}\,,\,\bX \,] := [\, \bV^{\lbrace \ell+1 \rbrace},...,\bV^{\lbrace d \rbrace}\,,\,\bX \,] - \bV^{\lbrace \ell \rbrace} \cdot [\bZ_\text{V}^{\lbrace \ell+1 \rbrace},...,\bZ_\text{V}^{\lbrace d \rbrace}\,,\,\bZ_\text{X}]$}
				\EndIf
			\EndIf
		\EndFor
		\EndProcedure
\end{algorithmic}

\subsection{Derivation of the algorithm}

\begin{figure}
	\centering
	\includegraphics[width=1\linewidth]{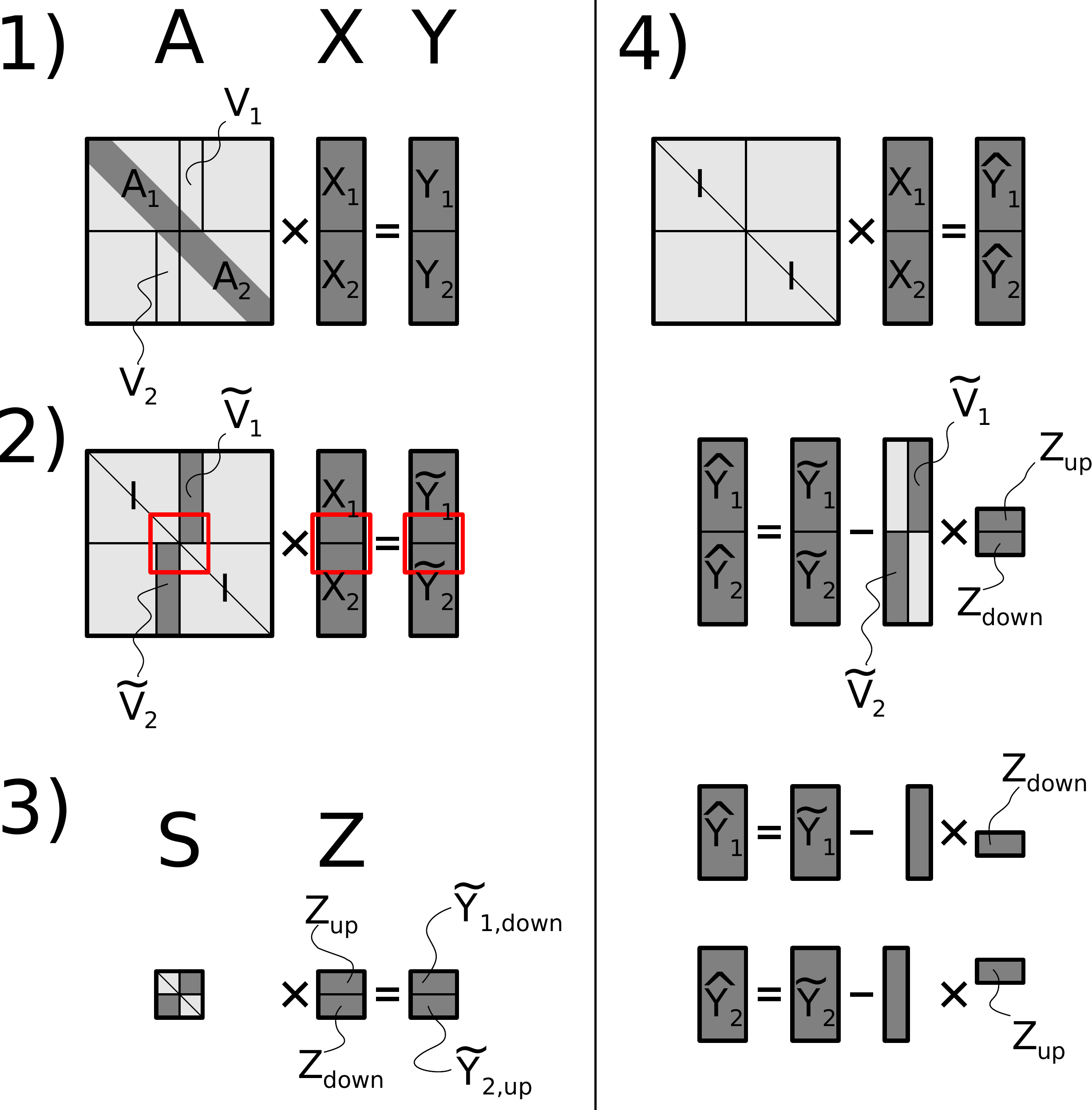}
	\caption{Derivation of our algorithm by applying SPIKE in a divide-and-conquer fashion. The system is divided into two equally sized systems of smaller dimension. After solving these, the solutions can be put together through the solution of a decoupled system (part 3) that has a small dimension.}
	\label{fig:Derivation_Spike}
\end{figure}

\paragraph{Origins in SPIKE}
We derived the above algorithm by applying the SPIKE algorithm \cite{10.1007/978-3-642-03869-3_74} due to Sameh in a divide-and-conquer fashion. We explain the derivation with the help of Figure~\ref{fig:Derivation_Spike}.

Initially, we are given a banded linear system $\bA\cdot \bX = \bY$, as is shown in part 1 of the figure. As is common for the divide-and-conquer approach, the system is split in the middle. We express the system as the composition of equally dimensioned square-matrices $\bA_1,\bA_2$, and matrices $\bX_1,\bX_2$, $\bY_1,\bY_2$. Since $\bA$ is banded, we need two additional matrices $\bV_1,\bV_2$ in order to be able to express $\bA$ as blocks. $\bV_1,\bV_2$ are high and thin: They have the hight of half the dimension of the original system and their breadth is equal to the bandwidth of the original system. We give names to $\tbV_1,\tbV_2$ the system in part 2: The system we call \textit{blade} and the matrices $\tbV_1,\tbV_2$ we call \textit{wings}. We will come back to this later.

In part 2 of the figure we see a transformed system that is obtained when multiplying the inverse of the block-diagonal matrix $\bD$
\begin{align*}
	\bD = \begin{bmatrix}
	\bA_1 & \bO\\
	\bO & \bA_2
	\end{bmatrix}
\end{align*}
from the left onto the system. If the dimension of $\bA_1$ and $\bA_2$ is large then it is unlikely that $\bD$ is singular. If $\bA$ is symmetric positive definite then $\bD$ is regular with a condition number bounded by that of $\bA$ \cite{Saad:2003:IMS:829576}. Thus, there exist conditions for $\bA$ such that the system considered in part 2 of the figure is well-posed. The obtained system has a interesting structure: It is an identity matrix plus two dense sub-blocks $\tbV_1,\tbV_2$, that have the breadth of $\bA$'s bandwidth.

In part 3 of the figure we consider a sub-system that is obtained when extracting the red portions from the matrices in part 2. This sub-system is decoupled from the total system. It has a size to twice the bandwidth and it can be solved directly (e.g. by $LU$-factorization followed by forward and backward substitutions) for the extracted (red-marked) portion of $\bX$. This approach would be followed in the usual SPIKE algorithm \cite{Polizzi2007}. However, we use a different approach. As shown in part 3 of our figure, we solve the sub-system with system matrix $\bS$ and write the solution into a matrix that we call $\bZ$. $\bZ$ is composed vertically of two blocks $\bZ_\text{up}$ and $\bZ_\text{down}$, of which each has the hight of the bandwidth of $\bA$, and the breadth of $\bX$ and $\bY$.

In part 4 we consider again a modified system. Using $\bZ$ from step 3, we see that we can change the right-hand sides such that the system matrix simplifies to an identity. This transformation can be easily derived: As the first block of the system in part 2, we have
\begin{align*}
	\bI \cdot \bX_1 + [\tbV_1,\bO] \cdot \bX_2 = \tbY_1\,.
\end{align*}
Replacing $[\tbV_1,\bO] \cdot \bX_2$ by $\tbV_1 \cdot \bZ_\text{down}$ and moving the second term to the right-hand side yields the formula for $\hbY_1$. The formula for $\hbY_2$ can be found in an analogous way.

\paragraph{Parallelism, recursions, and the computing system}
In Figure~\ref{fig:Derivation_Spike}, there are two stages where parallelism can be exploited: The transformation from part 1 to part 2 involves the solution of the following banded linear systems:
\begin{subequations}
\begin{align}
\bA_1 \cdot [\tbV_1,\tbY_1] &= [\bV_1,\bY_1]\\
\bA_2 \cdot [\tbV_2,\tbY_2] &= [\bV_2,\bY_2]
\end{align}\label{eqn:DC_linSys}
\end{subequations}
These two systems can be solved independently from each other. The second stage, where parallelism can be exploited, is in the computation of the two components $\hbY_1$ and $\hbY_2$:
\begin{subequations}
	\begin{align}
	\hbY_1 &= \tbV_1\cdot \bZ_\text{down}\\
	\hbY_2 &= \tbV_2\cdot \bZ_\text{up}
	\end{align}\label{eqn:DC_hY}
\end{subequations}

\begin{figure}
\centering
\includegraphics[width=0.7\linewidth]{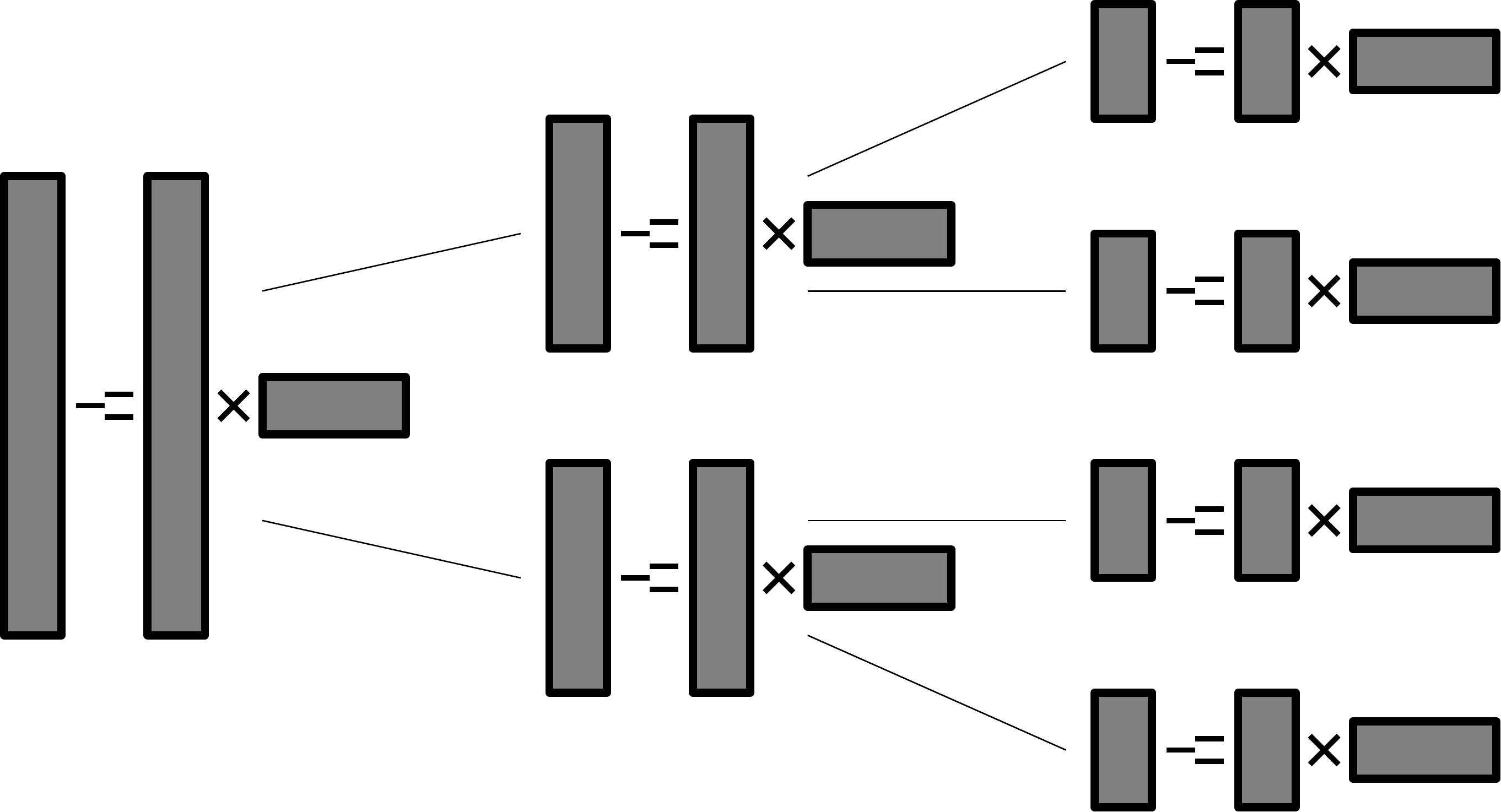}
\caption{Recursive divide-and-conquer approach to compute the matrix update $\bM -= \bU \cdot \bW$ on a parallel two-tree computing system. Each node divides the computational problem vertically into two problems of smaller dimension. These are then solved recursively by its children.}
\label{fig:Recursion_MatMat}
\end{figure}

The problems in \eqref{eqn:DC_linSys} and \eqref{eqn:DC_hY} can be expressed as recursions. The recursion for \eqref{eqn:DC_linSys} is obvious because we develop an algorithm that solves banded linear system and interiorly requires the solution of smaller banded linear systems. A recursion for \eqref{eqn:DC_hY} is found by distributing the computation over the rows of $\hbY_1,\hbY_2$. This is shown in Figure~\ref{fig:Recursion_MatMat}. The computation of each $\hbY_1$ and $\hbY_2$ can each be expressed in a form as $\bM := \bM - \bU \cdot \bW$, in some programming languages also written as $\bM -= \bU \cdot \bW$. The figure shows how the computation of the update of $\bM$ can be distributed through a two-tree network by dividing it in vertical direction.
All-together, we can formulate our whole algorithm for solving a banded linear system in a recursive way. The following code demonstrates this:
\begin{algorithmic}[1]
	\Procedure{RecursiveSolver}{$\bA,\bY$}
	\If($\bA$ has small dimension)
		\State $\bX = \bA \backslash \bY$
		\State \Return $\bX$
	\EndIf
	\State Decompose the system into $\bA_1,\bA_2,\bV_1,\bV_2,\bY_1,\bY_2,\bX_1,\bX_2$.
	\State $[\tbV_1,\tbY_1] := $\textsc{RecursiveSolver}($\bA_1,[\bV_1,\bY_1]$)
	\State $[\tbV_2,\tbY_2] := $\textsc{RecursiveSolver}($\bA_2,[\bV_2,\bY_2]$)
	\State Compose $\bS$ and solve the reduced linear system for $\bZ_\text{up}$, $\bZ_\text{down}$.
	\State \textit{// $\hbY_1$ \underline{is} $\tbY_1$, $\hbY_2$ \underline{is} $\tbY_2$}
	\State $\hbY_1 := $\textsc{RecursiveMatMul}($\tbY_1,\tbV_1,\bZ_\text{down}$)
	\State $\hbY_2 := $\textsc{RecursiveMatMul}($\tbY_2,\tbV_2,\bZ_\text{up}$)
	\State \textit{// $\bX$ \underline{is} $\begin{bmatrix}
			\hbY_1\\
			\hbY_2
		\end{bmatrix}$ }
	\State \Return $\bX$
	\EndProcedure
\end{algorithmic}

The computing system has been chosen as a two-tree in order to exploit the recursive nature of the algorithm: Initially, the root has the solve the linear system. According to the above code it would call its children in lines 7--8 to solve recursively the subsystems with $\bA_1$ and $\bA_2$. This is supposed to be done in parallel. Afterwards, the root composes $\bS$ from the small portions $\tbV_\text{1,down},\tbV_\text{2,up}$ of $\tbV_1,\tbV_2$, and computes $\bZ_\text{up},\bZ_\text{down}$. Finally, it calls again its children in order to compute $\hbY_1,\hbY_2$.

\paragraph{Access pattern on the distributed memory}
So far our explanation helps to understand how the algorithm works and why the results for $\bX$ are correct. But yet it is not easy to see how communication-intense the algorithm is and how the recursion acts on the global system. In this paragraph we give illustrations for both.

Our algorithm can be interpreted in an elegant way. To solve a system
\begin{align*}
	\underline{\bA} \cdot \underline{\bX} \cdot \underline{\bY}
\end{align*}
we can interprete that the algorithm applies an iterative scheme
\begin{align*}
	\bA^{(0)} &:= \underline{\bA} & \bY^{(0)} &:= \underline{\bY}\\
	\bA^{(j+1)} &:= \big(\bD^{(j)}\big)^{-1} \cdot \bA^{(j)} & \bY^{(j+1)} &:= \big(\bD^{(j)}\big)^{-1} \cdot \bY^{(j)}\quad \text{ for }j=0,...,d\,,
\end{align*}
where $d = \log_2(N)$ and where $\bA^{(d+1)} = \bI$ and thus $\underline{\bX} = \bY^{(d+1)}$. The matrices $\bD^{(j)}$ are block-diagonal matrices whose blocks are blades. As we have seen, the inverses of blades can be computed in parallel through \textsc{RecursiveMatMul}. We only need one upwards-communication of $\bV_\text{1,up}$, $\bV_\text{2,down}$ and one downwards communication of $\bZ_\text{down}$, $\bZ_\text{up}$.

\begin{figure}
\centering
\includegraphics[height=15cm]{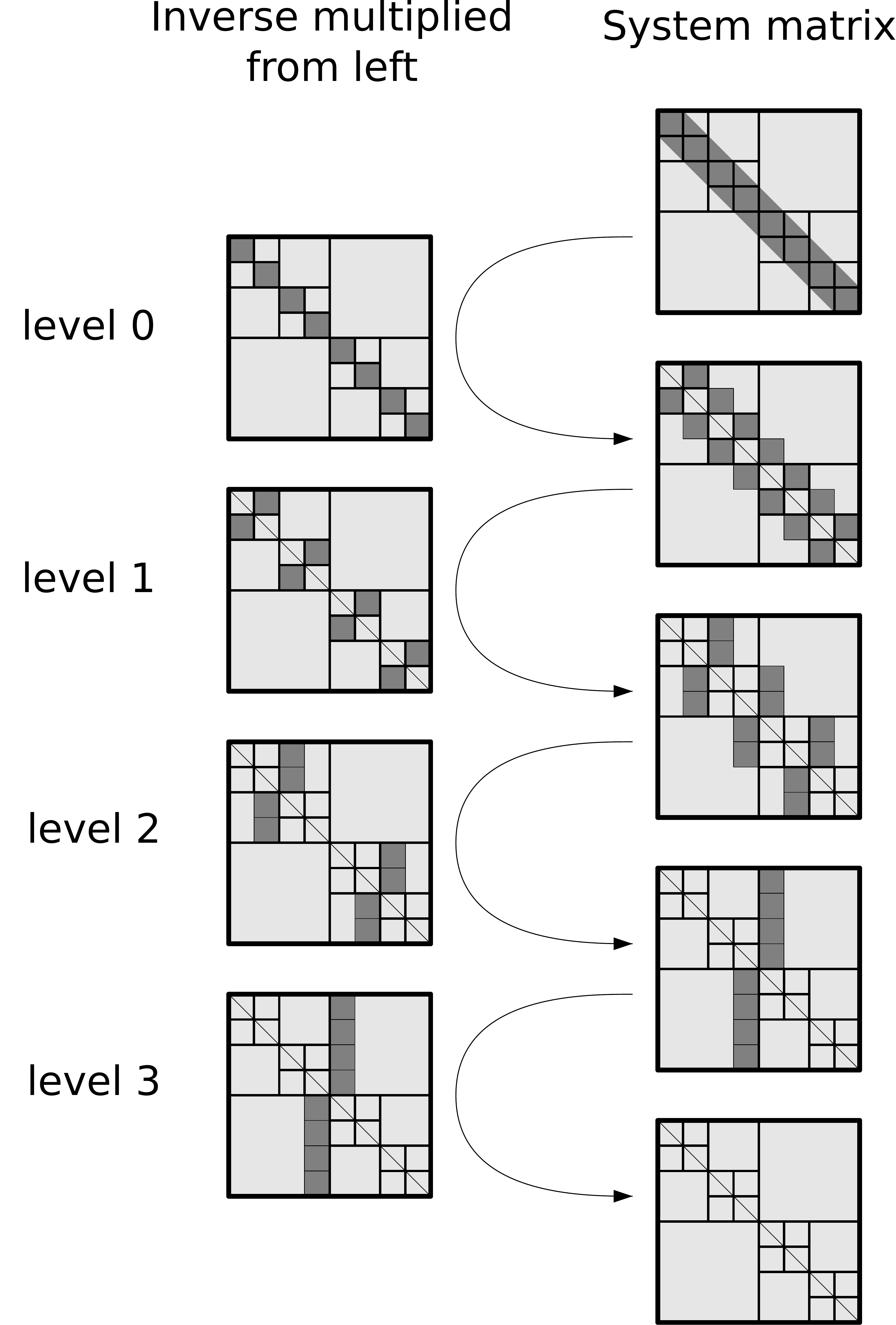}
\caption{Successive multiplication of inverse of block-diagonal matrix from the left onto the system matrix. The block-diagonal matrices consist of blades. Their inverse is easy to apply in parallel.}
\label{fig:Matrix_Pattern}
\end{figure}

Figure~\ref{fig:Matrix_Pattern} illustrates the iteration for $N=8$, $d=3$. The figure shows the system matrix in the right part, starting at the top with the original matrix and ending at the bottom with an identity. The left part of the figure shows the matrices $\bD^{(j)}$ for $j=0,...,d$\,.
\largeparbreak

We want to explain how this interpretation of the algorithm's action can be derived from \textsc{RecursiveSolver} and at the same time discuss the algorithmic steps as represented by the figure. The recursion in lines 7--8 results in a partitioning of the system matrix along the diagonal blocks. Black lines indicate the recursive diagonal partitioning. Grey coloring shows the non-zero pattern. In the bottom of the recursion the inverses of the smallest diagonal blocks are applied from the left onto the system. So the matrix $\bD^{(0)}$ consists solely of diagonal elements. The matrix $\bA^{(1)}$ has identity matrices on the diagonal blocks by construction. The sub-diagonal blocks are no longer triangular but dense. Since level 0 is the recursive bottom, we now ascend. The matrix $\bD^{(1)}$ consists of four blades. Multiplication from the left with its inverse yields a matrix $\bA^{(2)}$ that has identities on larger diagonal blocks (since these where identical to the blades on the diagonals of $\bD^{(1)}$). However, the hight of the wings in $\bA^{(2)}$ increases because there are subdiagonal blocks in $\bA^{(1)}$ that were not represented in $\bD^{(1)}$. Further ascending, the matrix $\bD^{(2)}$ consists of two diagonal blocks, which are blades of dimension $4$. The two subdiagonal blocks in $\bA^{(3)}$ suffer from fill-in in vertical direction while the diagonal blocks become identity matrices. Very finally, $\bA^{(3)}$ has a blade-structure and this $(\bD^{(3)})^{-1}\cdot\bA^{(3)}$ yields the identity. 

\largeparbreak

\begin{figure}
\centering
\includegraphics[width=1\linewidth]{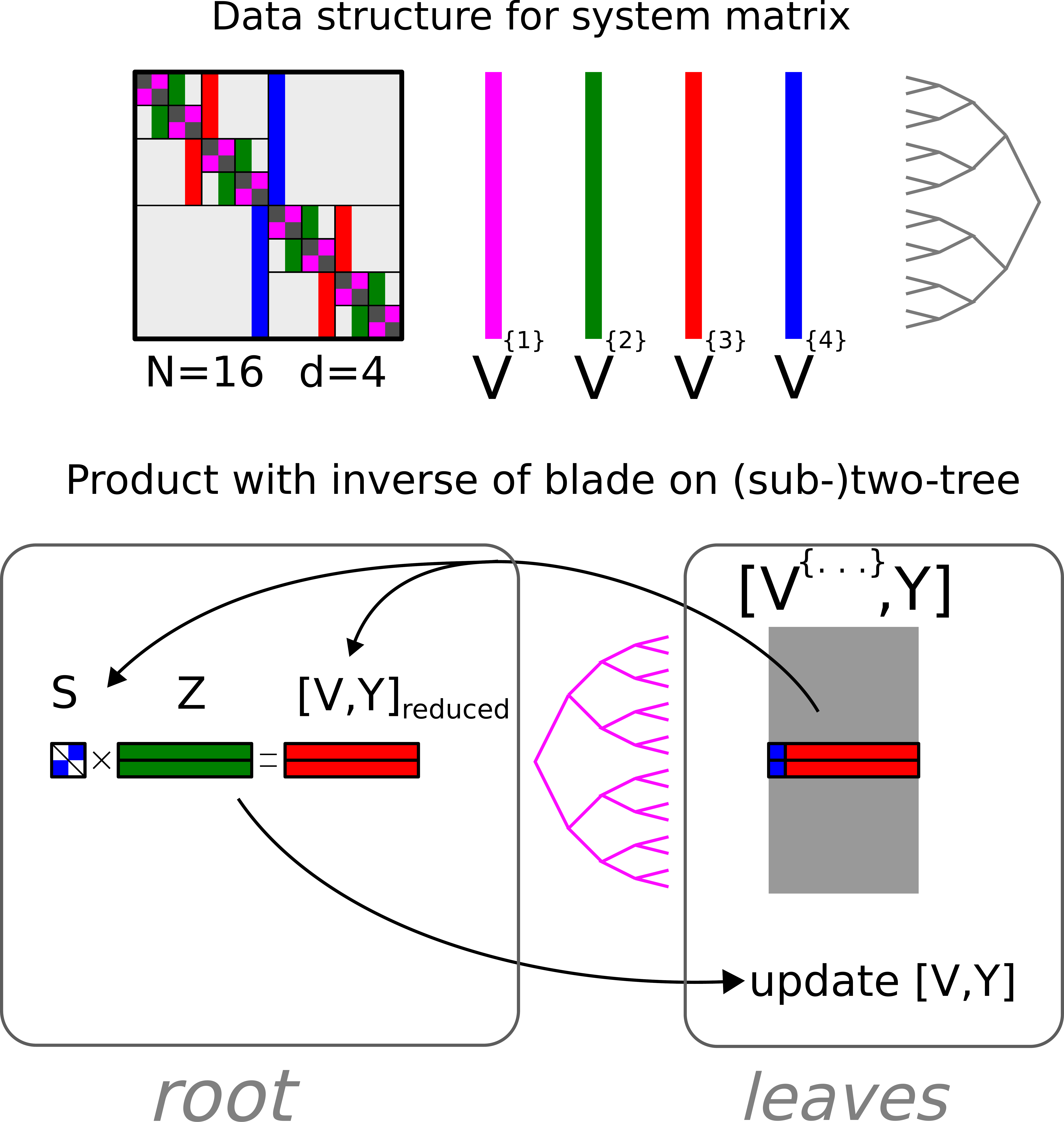}
\caption{Top: The matrix is represented by $d=\log_2(N)$ dense wing-matrices $\bV^{\lbrace j\rbrace}$, $j=1,...,d$, which are each distributed in row-chunks over the leaves. Bottom: The product with the inverse of a blade. The leaves send extracted matrices to the root of the sub-tree. The root computes the decoupled solution $\bZ$ and sends it back through the tree to all leaves, so they can update the wing-matrices.}
\label{fig:MemoryPattern}
\end{figure}

From the sparsity patterns of the system matrix in each iteration we can draw conclusions on the memory that the nodes need to hold in order to be able to compute $\bA^{(0)},...,\bA^{(d+1)}$ without the need of any overhead for, e.g., dynamically changing a sparse-memory representation of $\bA$. In our algorithm \textsc{ParallelSolver} we store $\underline{\bA}$ as wing-matrices $\bV^{\lbrace j}$, $j=1,...,d$\,. The top of Figure~\ref{fig:MemoryPattern} depicts this: The fill-in pattern of $\underline{\bA}$ fits into $\log_2(N)$ column vectors.

For this data structure the product of $\underline{\bA}$ with an inverse of on of the above block-diagonal matrices $\bD^{(j)}$ can be computed efficiently, as is shown in the bottom of the figure: Say we want to compute the product with $(\bD^{(j)})^{-1}$ for $j=3$. In this case, the leaves send the data for the reduced system (compare to the red-framed portions in Figure~\ref{fig:Derivation_Spike} part 2) to the root of the sub-trees of level $j=3$. The roots of the subtrees compute the reduced solutions $\bZ_\text{up},\bZ_\text{down}$, that afterwards they send back to all leaves through their sub-tree. The leaves update the wing-matrices by performing the same computations for them as for updating $\bY$.

\paragraph{Parallel time complexity}
We derive the time complexity by the following list of axioms:
\begin{enumerate}[label=(\roman*)]
	\item According to Figure~\ref{fig:Matrix_Pattern} the algorithm performs $\cO(\log(N))$ iterations.
	\item The communicated pieces of data per iteration are the node-local matrices
	\begin{align*}
		[\bV^{\lbrace 1 \rbrace},...,\bV^{\lbrace d \rbrace},\bY] \in \R^{ m \times (m \cdot d + k)}\,,
	\end{align*}
	which consist of $\cO\big( m \cdot (m \cdot \log(N)+k) \big)$ elements.
	The communication cost per iteration is bounded by the time required for sending these elements from an arbitrary leaf to the root or vice versa (size $\bZ_\text{V},\bZ_\text{X}$ have the same dimensions as $[\bV^{\lbrace 1 \rbrace},...,\bV^{\lbrace d \rbrace},\bY]$). Equivalently, the time complexity for communication per iterations is
	\begin{align*}
		\cO( c^\text{Lat}_N + m \cdot (m \cdot \log(N)+k ) )\,,
	\end{align*}
	where $c^\text{Lat}_N$ is the time complexity of physical time that is needed to send a single number from the root to a leaf that is farthest away from the root in terms of cable-length.
	\item The computational complexity per iteration per node is as follows: The nodes that are not leaves do either compute nothing or they compute a decomposition of a matrix $\bS \in \R^{(2 \cdot m) \times (2 \cdot m)}$ and apply it to the columns of the matrix $[\bV^{\lbrace 1 \rbrace},...,\bV^{\lbrace d \rbrace},\bY]$ in order to compute $\bZ_\text{V},\bZ_\text{X}$. The leaves on the other hand compute matrix-matrix-products of an $m \times m$-matrix with $\bZ_\text{V}$ and $\bZ_\text{X}$. Thus, the computational complexity per node per iteration is bounded by $\cO( m^3 \cdot \log(N) + m^2 \cdot k )$.
\end{enumerate}
Combining the items, we find that the parallel time complexity of our method is:
\begin{align}
	\cO\Bigg(\ \log(N)\cdot \Big(\,c^\text{Lat}_N+m^2 \cdot \big( m \cdot \log(N) + k \big) \,\Big)  \ \Bigg) \label{eqn:ComplexityOrder}
\end{align}
Unfortunately, for a computing system of $N$ nodes with each of a size in $\Theta(1)$ the minimum possible value for $c^\text{Lat}_N$ lives in $\Theta( N^{1/3} )$.

We consider two special cases:
\begin{enumerate}
	\item Assuming $k,m \in \cO(1)$ the complexity result simplifies to
	\begin{align*}
		\cO\Big(\ \log(N) \cdot N^{1/3} \ \Big)\,.
	\end{align*}
	\item Assuming $k,m \in \Omega(N^{1/6})$, the time complexity is
	\begin{align*}
		\cO\Big(\ \log(N) \cdot m^2 \cdot \big(m \cdot \log(N)+k\big) \ \Big)\,.
	\end{align*}
\end{enumerate}
Whereas the second result is obviously optimal for solving a band-matrix with dense band of bandwidth $m$ (since it is has $N$ diagonal blocks of size $m$ that need to be factorized), it turns out that the first result is not optimal and can be made optimal by a fine-tuning.

\section{Lower complexity bound for the parallel solution of tridiagonal linear systems}

In this section we prove in order the following statements, where $s$ is the physical dimension in which the computing system is built (e.g., if the computing system is built on the surface of a planet then $s=2$, and if the computing system is a/the planet then $s=3$):
\begin{enumerate}
	\item The latency $c^\text{Lat}_N$ is bounded from below by $\Omega(N^{1/s})$, regardless how many cables are used.
	\item The time complexity for an algorithm for solving the tridiagonal linear equation system of dimension $N$ on a distributed memory machine of computing nodes is bounded from below by $\Theta(N^{1/(s+1)})$, no matter how many nodes and cables are used.
\end{enumerate}

\paragraph{Latency}
We discuss on the latency time for a computing system of $N$ nodes that are $s$-dimensional spheres, where the latency of one node to communicate to another is bounded from below by the order of their physical distance in space.

For $N$ nodes, the diameter on a line is $\geq N-2$ because the two outermost nodes have $N-2 \in \Omega(N)$ nodes between themselves. For an illustration of this, consider Figure~\ref{fig:Block_Matrix_Storage}, that shows twice the topology for the two-tree computing system of $N=8$ nodes. Since the leaves are placed on a line with a unit distance, the physical distance of the root to the farthest leaf is bounded from below in $\Omega(N)$.

Now let us consider the case where we use $s=2$ dimensions: We place the computing nodes into the smallest possible circle. Since the total area of $N$ nodes is $\Theta(N)$, the radius $r$ of the must be in $r \in\Theta(N^{1/2})$. Then, the communication between two nodes requires at most a time of $2 \cdot r \in \Theta(N^{1/2})$. In $s=3$ dimensions the situation is even better. Here the radius must only be of length $r \in \Theta(N^{1/3})$. Given the positions of the nodes in the sphere, it is trivial to find an almost optimal communication network for them: Connecting them as a tree yields that the maximum diameter of the communication tree is $\cO(\,\log(N) \cdot N^{1/s})$ because each cable can at most have length $\cO(N^{1/s})$ and a tree network has diameter $\cO(\log(N))$.

\paragraph{Lower bound on time complexity for solving the tridiagonal linear systems}

Presume we shall solve a tridiagonal system $\bA \cdot \bx = \by$ of dimension $N$. Presume that we use $\cO(N^p)$ nodes of which each holds at most $\cO(N^q)$ pieces of data. It must be $q\geq 1-p$ because $\cO(N)$ pieces of data must be stored in total.

It is known that each number of the solution vector $\bx$ depends on each number of the right-hand side $\by$ and each value of the matrix $\bA$. Thus the following hold:
\begin{enumerate}[label=(\roman*)]
	\item Each node must read at least $\cO(N^{1-p})$ pieces of its data.
	\item Each node must communicate at least once (either directly or indirectly) to each other node.
\end{enumerate}
Combining the two properties, we find a lower bound for the time complexity with a fixed value of $p$:
\begin{align*}
	\Omega(\ \underbrace{N^{1-p}}_\text{read problem} + \underbrace{N^{\frac{p}{s}}}_\text{communicate} \ ) = \cO(N^\omega)
\end{align*}
The minimum of polynomial complexity orders $\omega$ depending on $s$ and the optimal order $p$ of the number of nodes are given in table \ref{table:ComplexityOrders}.

\begin{table}
	\centering
\begin{tabular}{||c|c|c||}
	\hline
	\hline dimension $s$ & order $p$ of number of nodes & order $\omega$ of time complexity \\
	\hline 	$1$ & $1/2$ & $1/2$ \\ 
	 		$2$ & $2/3$ & $1/3$ \\
	 		$3$ & $3/4$ & $1/4$ \\ 
	\hline 
	\hline
\end{tabular}
\caption[]{Lower complexity bounds for solving a linear equation system on parallel a distributed memory system of arbitrary many nodes.}\label{table:ComplexityOrders}
\end{table}

The question arises why we have not presented our algorithm with a number of nodes that is in $\cO(N^{3/4})$ since then our algorithm would be logarithmically close to the optimal time complexity. We did not because for our PhD thesis we will have to solve problems where $N \approx 10^6$ and $m \approx 100$. So for our applications the execution time is rather dominated by $\log(N)^2\cdot m^3$. The $m^3$ arises from the fact that $\underline{\bA}$ has $N$ dense $m \times m$-matrices on the diagonal that must be decomposed. In theory a fast matrix-multiplication algorithm could be employed to reduce the complexity order for this, but this is not practical. Thus, for the problems that we need to solve, the complexity of our algorithm is already logarithmically close to optimal or maybe even optimal.

\section{Conclusions}
We presented an algorithm for the efficient parallel solution of (block-)tridiagonal linear systems. We provided an accurate implementation of the algorithm that uses a common message-passing interface. Following our analysis, the algorithm has a parallel time complexity that could be made optimal for tridiagonal systems (by simply using fewer nodes) and that is clearly optimal for block-tridiagonal linear systems, respectively banded linear systems of small bandwidth.

Though a proper implementation has been provided for the algorithm, it will be rather difficult to utilize its full potential in practice. This is because the time that is required to send the data to the solver would already destroy the benefit. Software that uses this solver needs to be highly sophisticated. In particular, problems must be instantiated in a way such that the linear system is already distributed in the memory of the leaves of our computing system at the time when our solver is called.

Further work will be related to an attempt of implementing an optimal control solver by direct transcription that shall solve the linear systems within the non-linear programming solver by means of the linear system solver that has been presented in this work.

\FloatBarrier

\bibliography{recursiveSpikeBib}
\bibliographystyle{plain}

\end{document}